\documentclass{article}
\usepackage[utf8]{inputenc}
\usepackage[T1]{fontenc}
\usepackage{lmodern}
\usepackage{babel}
\usepackage{lipsum}
\usepackage{graphicx}
\usepackage{amsmath}
\usepackage{amsfonts}
\usepackage{subcaption}
\usepackage{tabularx}
\usepackage{siunitx}
\usepackage{booktabs}
\usepackage{tikz}
\usetikzlibrary{calc}
\usetikzlibrary{patterns}
\usepackage[ruled,linesnumbered]{algorithm2e}
\usepackage{accents}
\usepackage{parskip}
\usepackage{amssymb}
\usepackage[hyperfootnotes=false]{hyperref}
\usepackage{multirow}
\usepackage{alphabeta} 
\usepackage{todonotes}

\hypersetup{
    colorlinks=true,
    linkcolor=blue,
    filecolor=magenta,      
    urlcolor=cyan,
    citecolor=blue
}

\newcommand{\coefficients}{\ensuremath{\boldsymbol{\alpha}}\xspace}

\newcommand{\rhs}{\ensuremath{\beta}\xspace}

\newcommand{\xv}{\ensuremath{\mathbf{x}}\xspace}
\newcommand{\lpoptimal}{\ensuremath{\mathbf{x}^{LP}}\xspace}
\newcommand{\lpoptimali}[1]{\ensuremath{x^{LP}_{#1}}\xspace}

\newcommand{\reals}{\ensuremath{\mathbb{R}}\xspace}

\newcommand{\splitdisjunction}[2]{\ensuremath{\mathcal{D}(\boldsymbol{\pi},\pi_{0})}\xspace}

\newcommand{\psc}[3]{\ensuremath{\mathtt{psc}(#1,#2,#3)}\xspace}
\newcommand{\loc}[1]{\ensuremath{\mathtt{loc}(#1)}\xspace}
\newcommand{\sps}[1]{\ensuremath{\mathtt{sps}(#1)}\xspace}

\newcommand{\psvar}[1]{\ensuremath{\boldsymbol{\psi}(x_{#1})}\xspace}
\newcommand{\downlocksi}[1]{\ensuremath{\boldsymbol{\zeta}^{-}_{#1}}\xspace}
\newcommand{\uplocksi}[1]{\ensuremath{\boldsymbol{\zeta}^{+}_{#1}}\xspace}

\newif\ifarxiv
\arxivtrue

\ifarxiv
\usepackage[scale=0.9]{FiraMono}
\usepackage{arxiv}
\usepackage{charter}
\usepackage{fullpage}
\title{A Context-Aware Cutting Plane Selection Algorithm for Mixed-Integer Programming}
\author{ \href{https://orcid.org/0000-0001-7270-1496}{\includegraphics[scale=0.06]{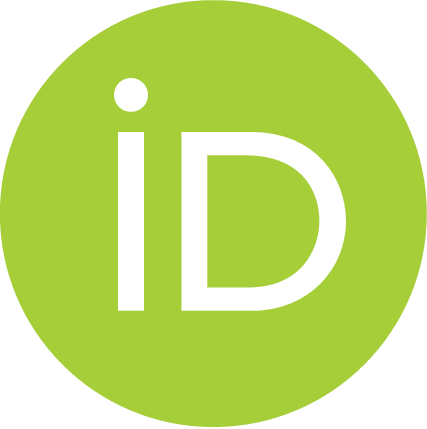}}\hspace{1mm}Mark Turner\thanks{Institute of Mathematics, Technische Universit{\"a}t Berlin, Straße des 17. Juni 135, 10623 Berlin, Germany}\hspace{2mm}\thanks{Zuse Institute Berlin, Department of Mathematical Optimization, Takustr. 7, 14195 Berlin} \\
\texttt{turner@zib.de} \\
\And
\href{https://orcid.org/0000-0002-6320-8154}{\includegraphics[scale=0.06]{orcid_id_icon.eps}}\hspace{1mm}Timo Berthold\footnotemark[1]\hspace{2mm}\thanks{Fair Isaac Germany GmbH, Takustr. 7, 14195 Berlin, Germany} \\
\texttt{timoberthold@fico.com}
\And
\href{https://orcid.org/0000-0002-6284-3033}{\includegraphics[scale=0.06]{orcid_id_icon.eps}}\hspace{1mm}Mathieu Besançon\footnotemark[2]\\
\texttt{besancon@zib.de} \\
}
\else
\title{A Context-Aware Cutting Plane Selection Algorithm for Mixed-Integer Programming}
\author{Mark Turner$^*$\inst{1,2}\orcidID{0000-0001-7270-1496} \and
Timo Berthold\inst{1,3}\orcidID{0000-0002-6320-8154}\and\\
Mathieu Besançon\inst{2}\orcidID{0000-0002-6284-3033}}
\authorrunning{M.~Turner, T.~Berthold, M.~Besançon}
\institute{Institute of Mathematics, Technische Universit{\"a}t Berlin, Straße des 17. Juni 135, 10623 Berlin, Germany \and
Zuse Institute Berlin, Department of Mathematical Optimization, Takustr. 7, 14195 Berlin \and
Fair Isaac Germany GmbH, Takustr. 7, 14195 Berlin, Germany \\
\email{turner@zib.de}, \email{timoberthold@fico.com}, \email{besancon@zib.de}}
\fi

\begin{document}

\maketitle

\begin{abstract}
The current cut selection algorithm used in mixed-integer programming solvers has remained largely unchanged since its creation.
In this paper, we propose a set of new cut scoring measures, cut filtering techniques, and stopping criteria, extending the current state-of-the-art algorithm and obtaining a 5\% performance improvement for SCIP over the MIPLIB 2017 benchmark set.
\ifarxiv
\else
\fi
\end{abstract}

\section{Introduction and Related Work}

\renewcommand*{\thefootnote}{\fnsymbol{footnote}}
\footnotetext[1]{Corresponding author.}
\renewcommand*{\thefootnote}{\arabic{footnote}}

We focus on this paper on the selection of cutting planes, an essential component of the branch-and-cut framework, which is the main algorithmic paradigm to solve Mixed-Integer Linear Programs (MILP) classically defined as:
\begin{align}
    \underset{\mathbf{x}}{\mathrm{argmin}}\{\mathbf{c}^{\intercal}\mathbf{x} \;\; | \;\; \mathbf{A}\mathbf{x} \leq \mathbf{b}, \;\; \mathbf{l} \leq \mathbf{x} \leq \mathbf{u}, \;\; \mathbf{x} \in \mathbb{Z}^{|\mathcal{J}|} \times \mathbb{R}^{n - |\mathcal{J}|} \} \label{eq:mip}
\end{align}
A cut, parameterised by $(\coefficients, \rhs) \in \reals^{n+1}$, is an inequality $\coefficients^{\intercal} \xv \leq \beta$ that is violated by at least one solution of the LP relaxation but that does not increase the optimal value of the problem when added i.e.~it is valid for \eqref{eq:mip}.
It is used to tighten the Linear Programming (LP) relaxation of \eqref{eq:mip}.
Cuts are generated in rounds, where we refer to the process of computing cuts, applying a subset of the computed cuts, and resolving the LP as a \emph{separation round}. In general, cuts are cheap to compute, and more are computed than actually applied to the LP. The process of deciding which cuts to add to the LP is called \emph{cut selection}, and is the focus of this paper, in which we introduce a new cut selection algorithm. 

Cut selection was considered largely ``solved", following computational results on a variety of MILP solvers \cite{achterberg2007constraint,zerohalf,wesselmann2012implementing}, in that cheap heuristic selection rules were sufficient for good performance. Recently this conclusion has been challenged, with largely machine learning-driven research attempting to extract further performance, see \cite{deza2023machine} for an overview.
Such research, however, is often limited in how it can be deployed in a MILP solver due to the complexity of deploying black-box predictors and their brittle generalisation.
In this paper we introduce a new cut selection algorithm, which revises the three major aspects of cut selection, namely cut scoring measures, cut filtering techniques, and stopping criteria.
We determine default parameter values of our algorithm by leveraging methodology \cite{hutter2010sequential,hutter2010automated,xu2011hydra} and software \cite{smac} from the algorithm configuration community.
Our algorithm and the accompanying default values will be integrated in the next release of SCIP \cite{scip8}. 

\section{New Cut Selection Techniques}

Our new cut selection algorithm expands and generalises the existing one from SCIP \cite{scip8} in three aspects, namely, cut scoring, cut filtering, and stopping criteria.
For all three aspects, we incorporate additional information from the cut and from the current context in which the separation takes place.

\subsection{Cut Scoring} \label{subsec:cut_scoring}

The most studied aspect of cut selection is the scoring of cuts, which subsumes their ranking. This is both due to the large variety of scoring measures that exist, see \cite{wesselmann2012implementing,turner2023cutting}, and the common assumption of a fixed selection rule outside of scoring \cite{tang2020reinforcement,balcan2022improved,turner2022adaptive}. We introduce a new set of cut scoring measures that incorporate information from other non-cutting plane aspects of the MILP solving process.

\subsubsection{Pseudo-cost Scoring}

Pseudo-costs \cite{achterberg2005branching}, the dominant decision-maker in state-of-the-art MILP branching rules \cite{achterberg2009hybrid}, estimate scores for branching candidates based on the historical objective value improvement observed when branching on the candidates. It has recently been shown in \cite{turner2023branching} that cut scoring measures can complement pseudo-cost-based branching rules for improved solver performance. We address the reverse direction here, and use pseudo-costs to augment cut scores, where we let \psvar{i} be the pseudo-cost of variable $x_{i}$.
We define the pseudo-cost scoring of a cut, $(\coefficients, \rhs) \in \reals^{n+1}$, as the expected objective change as predicted by the pseudo-costs and express it as follows:
\begin{align}
    \psc{\coefficients}{\rhs}{\lpoptimal} := \sum\limits_{{i, \alpha_{i} \neq 0}} \psvar{i} \left(|\lpoptimali{i} - \alpha_{i} \frac{\coefficients^{\intercal} \lpoptimal  - \beta}{\|\coefficients\|^{2}}|\right)
\end{align}
We normalise the result by the maximum pseudo-cost score of a cut in the separation round, so all pseudo-cost scores are in the range $[0,1]$. 

\subsubsection{Lock Scoring}

Down-locks and up-locks, introduced in \cite{achterberg2007constraint}, can be interpreted for a variable $x_{i}$, as the number of constraints that ``block" the shifting of $x_{i}$ in the direction of the variable's lower or upper bound. Let \downlocksi{i} be the number of down-locks and \uplocksi{i} be the number of up-locks for variable $x_{i}$. We define the lock scoring of a cut, $(\coefficients, \rhs) \in \reals^{n+1}$, as:
\begin{align}
    \loc{\coefficients} :=  \frac1n \sum\limits_{i, \alpha_i \neq 0}({ \downlocksi{i} + \uplocksi{i}})
\end{align}
As with the pseudo-cost measure, we normalise the result by the maximum among all cuts in the separation round so the scores are in the range $[0,1]$. Additionally, as it is initially unclear if we want to reward a cut featuring variables with many locks, i.e.~further restrict heavily restricted variables, we also introduce the complement of this measure, i.e.~promote cuts that restrict variables with few restrictions.

\subsubsection{Sparsity Scoring}
Sparsity refers to the fraction of zero coefficients in a cut. In general, dense cuts slow down LP solves \cite{sparsity}. We therefore introduce a score that promotes sparse cuts, which, given a cut $(\coefficients, \rhs) \in \reals^{n+1}$, and the parameters $\texttt{maxsps} \in \reals_{+}$, and $\texttt{endsps} \in [0,1]$, is defined as:
\begin{align}
    \sps{\coefficients} := \max \{\texttt{maxsps} -\frac{\texttt{maxsps}}{\texttt{endsps}} \cdot \; \| \alpha \|_0 , \;\; 0 \} 
\end{align}
\subsection{Cut Filtering} \label{subsec:cut_filtering}

The standard approach in the literature for cut filtering is the parallelism-based approach \cite{achterberg2007constraint,zerohalf,wesselmann2012implementing}. The approach selects the highest-scoring cut, removes any remaining cut from consideration that is too parallel to the selected cut, and then repeats until some stopping criteria is met or no more non-filtered cuts are left to be selected. Such a procedure was shown to be absolutely necessary for good performance in \cite{achterberg2007constraint,zerohalf}, however that is only w.r.t.~trivial selection approaches. We therefore present and discuss two other approaches for cut filtering. 

\subsubsection{Density Filtering}

Density-based filtering, unlike parallelism-based filtering, is performed at the start of the selection process. It immediately removes any generated cut above a given density threshold. Preliminary results for such an approach already exist, see \cite{turner2023cutting,herman}, and suggest that $40\%$ density, i.e.~a cut with at most $40\%$ of its entries being non-zero, is a potentially desirable setting.

\subsubsection{Parallelism-Based Penalties}

It is possible to simply impose a penalty on any cut parallel to previously-added ones as opposed to outright removing it from consideration. In our algorithm, this involves reducing the score of all remaining cuts considered parallel to a selected cut.
If the score falls below some threshold, which we set to 0, then the cut is removed from consideration. 

\subsection{Stopping Criteria}

In general for computational studies, and MILP solvers by extension, there is a hard-coded limit on the number of cuts that can be added per round. Typically, this limit is not reached, with the standard parallelism filtering algorithm removing enough cuts by itself as detailed in Subsection~\ref{subsec:cut_filtering}. As we are employing alternate filtering algorithms, however, we believe it is necessary to have an instance-dependent limit. We therefore place a limit on the number of \emph{non-zeros} added to the LP at each separation round, where the limit is a multiple of the number of columns in the LP.

\section{Experiments and Computational Results}\label{sec:experiments}

We conduct the following two experiments: First we select an appropriate training set for determining the best algorithm configuration of our new cut selector. Second we deploy the best found configuration on a test set to determine the actual MILP solver performance improvement. 

Our training set is a combination of instances from the MIPLIB 2017 collection set\footnote{MIPLIB 2017 -- The Mixed Integer Programming Library \url{https://miplib.zib.de/}.} \cite{miplib}, strIPlib\footnote{\url{https://striplib.or.rwth-aachen.de/login/}} \cite{striplib}, and SNDlib-MIPs \cite{zenodo_sndlib}. We select the training set by removing any instance from consideration that takes less than 5 seconds or longer than 120 seconds to solve, takes longer than 10 seconds to presolve, or solves in less than 50 or more than 20000 nodes. This selection procedure is repeated with an optimal primal solution preloaded. From the remaining instances, we create a mapping to a feature space, where the features are the proportion of each constraint and variable type, the average row density, the maximum row density, and the objective vector density. Using this mapping, we select 40 instances that are maximally diverse by maximising the distance between any two selected instances in the feature space. Our test set is the MIPLIB 2017 benchmark set.

For all experiments, SCIP 8.0.3 \cite{scip8} is used, with random seeds $\{1,2,3,4,5\}$, PySCIPOpt \cite{pyscipopt} as the API, and Xpress 9.0.2 \cite{xpress} as the LP solver restricted to a single thread. For training, we use in non-exclusive mode a cluster equipped with Intel Xeon Gold 6342 CPUs running at 2.80GHz. For testing we run in exclusive mode on a cluster equipped with Intel Xeon Gold 5122 CPUs running at 3.60GHz, where each run is restricted to 48GB memory and a 2-hour time limit. The code used for all experiments is available and open-source\footnote{\url{https://github.com/Opt-Mucca/Context-Aware-Cut-Selection}}, and will be integrated in the next release of SCIP.

In addition to the new cut scoring measures introduced in Subsection~\ref{subsec:cut_scoring}, we use standard scoring measures from the literature. Namely, we use efficacy, $\texttt{eff}$, expected improvement, $\texttt{exp}$, objective parallelism, $\texttt{obp}$, and integer support, $\texttt{isp}$ \cite{wesselmann2012implementing}. We normalise efficacy and expected improvement as in \cite{turner2022adaptive}, ensuring values in the range $[0,1]$. Notably, we disable directed cutoff distance due to it being unfairly advantaged by preloading the optimal solution in training.
All cut scoring measures are included in a weighted sum rule, where the multipliers are hyperparameters. 

\subsection{Training}

We use SMAC \cite{smac} to determine the parameter values of our new cut selection algorithm. We run 5 instances of SMAC on our training set using the random seeds 1-5 and a limit of 300 complete passes of SMAC over our training set. Our training metric is the ratio of shifted geometric means of solve time compared to SCIP default with a shift of one second. 

\begin{table}[h]
\centering
    \begin{tabular}{l|ccccc}
        \hline
        \hline
        Seed & 1 & 2 & 3 & 4 & 5 \\
        \hline
        \hline
        \noalign{\vskip 0.35mm}
        Time (s) & 0.874 & 0.909 & 0.853 & 0.917 & 0.891 \\
        Nodes & 0.847 & 0.951 & 0.862 & 0.920 & 0.910 \\
        \hline
        \hline
    \end{tabular}
    \caption{Ratio of shifted geometric means (1s for solve time and 10 for nodes) for our approach vs SCIP default on instance-seed pairs of training set.}
    \label{tab:head_to_head_geom_ratio_smac}
\end{table}
\ifarxiv
\else
\vspace{-1cm}
\fi
The runs generated largely different parameter choices, with all choices resulting in an 8-15\% improvement w.r.t. solve time over SCIP default on the training set. The generated parameter choices differed mostly in how to score a cut, where we could find no common parameter choices. For filtering and penalising cuts, however, there was an overwhelming consensus for the following parameters: dense cuts should be filtered at some $40-50\%$ threshold, and parallel cuts should not be filtered, but rather lightly penalised. We additionally note that due to the instance restrictions put on our training set, we do not expect the level of performance to perfectly generalise to larger MILP instances. 

\subsection{Test}

As the cut selection performance space is highly non-linear, aggregating the parameter choices from the various SMAC runs is difficult. We therefore use a combination of the results from SMAC (e.g.~density and parallelism filtering), and expert knowledge on the impact of cuts on the overarching solving process. 

\begin{table}[h]
    \centering
    \begin{tabular}{ccccc}
        \hline
        \hline
        \noalign{\vskip 1mm}
        \multicolumn{5}{c}{MIPLIB 2017 benchmark set} \\
        \multicolumn{5}{c}{97\% instance-seed pairs affected} \\
        \noalign{\vskip 1mm}
        \hline
        \noalign{\vskip 0.35mm}
        Time (s) & Nodes & $\Delta-$Solved & $\Delta-$Dual & $\Delta-$Primal \\
        \hline
        \textbf{0.95} & \textbf{0.92} & \textbf{14(/668)} & \textbf{1(/532)} & -11(/532) \\
        \hline
        \hline
    \end{tabular}
    \caption{Ratio of shifted geometric means (10s for solve time and 100 for nodes) for our approach vs SCIP default over affected instance-seed pairs of the test set (solved to optimality for both cut selection algorithms). $\Delta$ is the difference in wins over the instance-seed pairs for amount solved, and the respective bounds for unsolved instances. Bold entries are when better than SCIP default.}
    \label{tab:head_to_head_geom_ratio_ensemble}
\end{table}

\begin{figure}[h]
    \centering
    \includegraphics[width=0.75\textwidth]{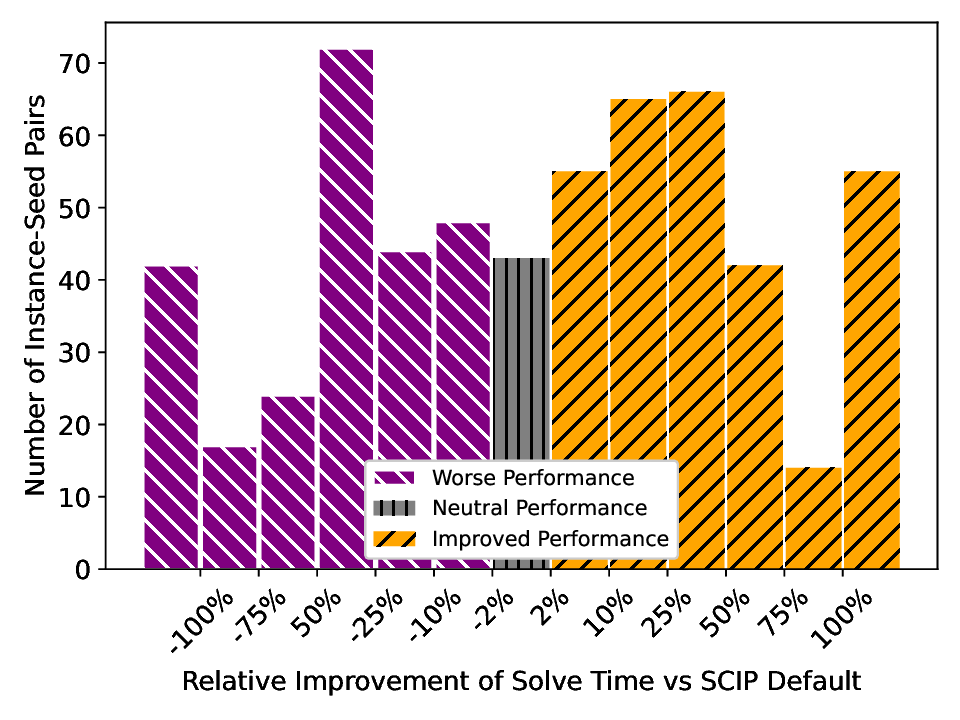}
    \caption{Bar plot of relative improvement of our approach vs SCIP default over affected instance-seed pairs of MIPLIB 2017 benchmark set.}
    \label{fig:bar_plot_ensemble}
\end{figure}

Comparing our new cut selector to SCIP default, we observe in Table~\ref{tab:head_to_head_geom_ratio_ensemble} an improvement of $5\%$ for solve time and $8\%$ for number of nodes.
The distribution of relative time improvement over the instance-seed pairs is visualised in Figure~\ref{fig:bar_plot_ensemble}, and is verified as statistically significant by a Wilcoxon signed-rank test (with a p-value  below 5\%).

\section{Conclusion}

In this paper we introduced a new cut selection algorithm, which combines new context-aware scoring measures, filtering methods, and stopping criteria. We incorporate the algorithm into SCIP, and obtain a 5\% improvement in solve time over the MIPLIB 2017 benchmark set.

\section*{Acknowledgements}
We thank Antonia Chmiela and Michael Winkler for helpful discussions about the experimental design. The work for this article has been conducted in the Research Campus MODAL funded by the German Federal Ministry of Education and Research (BMBF) (fund numbers 05M14ZAM, 05M20ZBM). The described research activities are funded by the Federal Ministry for Economic Affairs and Energy within the project UNSEEN (ID: 03EI1004-C).

\ifarxiv
\bibliographystyle{unsrt}
\bibliography{mybib}
\else
\bibliographystyle{splncs04}
\bibliography{mybib}
\fi

\end{document}